\documentclass[11pt]{article}
\usepackage{amssymb,amsfonts,amsmath,latexsym,epsf,tikz,url}
\newtheorem{theorem}{Theorem}[section]
\newtheorem{proposition}[theorem]{Proposition}
\newtheorem{observation}[theorem]{Observation}

\newtheorem{corollary}[theorem]{Corollary}
\newtheorem{lemma}[theorem]{Lemma}

\newtheorem{example}[theorem]{Example}

\newcommand{\proof}{\noindent{\bf Proof.\ }}
\newcommand{\qed}{\hfill $\square$\medskip}
\allowdisplaybreaks

\begin{document}

\title{On the saturation number of graphs}

\author{Saeid Alikhani$^{}$\footnote{Corresponding author} and Neda Soltani}

\date{\today}

\maketitle

\begin{center}

Department of Mathematics, Yazd University, 89195-741, Yazd, Iran\\
{\tt alikhani@yazd.ac.ir, neda\_soltani@ymail.com }\\

\end{center}


\begin{abstract}
Let $G=(V,E)$ be a simple connected graph. A matching $M$ in a graph $G$ is a collection of edges of $G$ such that no 
two edges from $M$ share a vertex. A matching $M$ is maximal if it cannot be extended to a larger matching in $G$. The cardinality of any smallest maximal matching in $G$ is the saturation number of $G$ and is denoted by $s(G)$. 
In this paper we study the saturation number of the corona product of two specific graphs. We also consider some 
graphs with certain constructions that are of importance in chemistry and study their saturation number. 
   
\end{abstract}


\medskip
\noindent{\bf AMS Subj.\ Class.:} 05C30; 05C70.


\section{Introduction}

Matching theory is a branch of graph theory concerned with study of structural and enumerative aspects of matchings, collections of edges of a graph that do not share a vertex. Its development has been strongly influenced and stimulated by chemical applications, in particular by the study of perfect matchings in benzenoid graphs. Additional impetus came with discovery of fullerenes, again mostly dealing with perfect matchings \cite{Ful, Sli, Kar, Zha}, but including also some structural results \cite{Ves, Sat}. For a general background on matching theory we refer the reader to the classical monograph by Lov\'asz and Plummer \cite{Lov}; for graph theory terms not defined here we also recommend \cite{Wes}.

The cardinality of matching $M$ of a graph $G$ is called the size of the matching. It is intuitively clear that matchings of small size are not interesting (each edge is a matching of size one, and the empty set is a matching of size $0$). Hence, we will be interested in matchings that are as large as possible. The cardinality of any maximum matching in $G$ is denoted by 
$\alpha^{\prime}(G)$ and is called the matching number of $G$. Since each vertex can be incident to at most one edge of 
a matching, it follows that the matching number of a graph on $n$ vertices cannot exceed $\left\lfloor\frac{n}{2}\right\rfloor$. If each vertex of $G$ is incident with an edge of $M$, the matching $M$ is called perfect. So the number of vertices of a graph $G$ admitting a perfect matching is even, but the opposite is generally not true.
Perfect matchings are obviously also maximum matchings. The study of perfect matchings, also known as Kekul\'e structures has a long history in both mathematical and chemical literature. For more details on perfect matching, we refer the reader to see \cite{Lov}. 

There is, however, another way to quantify the idea of large matchings. A matching $M$ is maximal if it cannot be extended to a larger matching in $G$. Obviously, every maximum matching is also maximal, but the opposite is generally not true.

Maximal matchings are much less researched than their maximum counterparts. That goes both for their structural and their enumerative aspects. While there is vast literature on perfect and maximum matchings (see, for example, monographs \cite{Lov} and \cite{Cyv}), the results about maximal matchings are few and scattered through the literature. We mention here two papers that treat, among other topics, maximal matchings in trees \cite{Kla, Wan}, one concerned with the structure of equimatchable graphs \cite{Fre}, and a recent paper about the saturation number of benzenoid graphs \cite{Dos}.


\medskip
The cardinality of any smallest maximal matching in $G$ is the saturation number of $G$ and is denoted by $s(G)$ (the same term, saturation number, is also used in the literature with a different meaning; we refer the reader to \cite{Fad} for more information). It is easy to see that the saturation number of a graph $G$ is at least one half of the matching number of $G$, i.e., $s(G)\geq \frac{\alpha^{\prime}(G)}{2}$ (\cite{Dos}). We recall that a set of vertices $I$ is independent if no two vertices from $I$ are adjacent. Clearly, the set of vertices that is not covered by a maximal matching is independent \cite{Edm}. This observation provides an obvious lower bound on saturation number of the graph $G$, i.e. $s(G)\geq \frac{(n-|I|)}{2}$, where $G$ is graph of order $n$ (\cite{Ves}).

\medskip 
Next section is concerned with the saturation number of some specific graphs. we also study the saturation number of the corona product of two certain graphs. In Section 3, we consider some graphs with specific construction that are of importance in chemistry and study their saturation number. As consequences, we obtain the saturation number of some cactus chains. 

\section{Saturation number of some graphs}

This section is organized as follows. In the Subsection 2.1, we consider the saturation number of some certain graphs.
We also examine the effects on the saturation number of them when they are modified by operations on their edge. Subsection 2.2 moves on to the saturation number of the corona product of two specific graphs.  

\subsection{Saturation number of some specific graphs} 

We start this section with the saturation number of some certain graphs. Let $P_n$, $C_n$ and $W_n$ be a path, cycle 
and wheel of order $n$, respectively. We have 
     \begin{enumerate}
            \item[(i)] $s(P_n)=\left\{
\begin{array}{lcl}
\left\lceil{\frac{n}{3}}\right\rceil, & \quad\mbox{if $n\equiv 2 ~ (mod\, 3)$}\\[15pt]
\left\lfloor{\frac{n}{3}}\right\rfloor, & \quad\mbox{otherwise.}
\end{array}
\right.$
           \item[(ii)] $s(C_n)=\left\lceil{\frac{n}{3}}\right\rceil$.
           \item[(iii)] $s(W_n)=1+s(P_{n-2})$.
      \end{enumerate}



 

           

\medskip
The union $G_1\cup G_2$ of graphs $G_1$ and $G_2$ with disjoint point sets $V(G_1)$ and $V(G_2)$ and edge sets $E(G_1)$ and $E(G_2)$ is the graph with $V(G_1)\cup V(G_2)$ and $E(G_1)\cup E(G_2)$ \cite{Har}. In addition, a maximum matching in a disconnected graph consists of the union of maximum matchings in each of its components. So we can 
conclude the following usefull lemma.

\begin{lemma}\label{union}
For any two simple graphs  $G_1$ and $G_2$,   $s(G_1\cup G_2)=s(G_1)+s(G_2).$
\end{lemma}

\medskip
Here we consider the effect on the saturation number of the path and the cycle graphs when one edge of them is deleted.

\begin{proposition}
    \begin{enumerate}
           \item [(i)] For every $n\geq 3$ and $e\in E(G)$, we have $$s({C_n}-e)=s(P_n)\leq s(C_n).$$
           \item [(ii)] Let $V(P_n)=\{v_1,v_2,...,v_n\}$ and $e_i=v_iv_{i+1} \in E(P_n)$ where $1\leq i\leq n-1$. For every $n\geq 1$, we have 
\begin{equation*}  
s({P_n}-{e_i})=\left\{
\begin{array}{ll}   
s(P_n)-1, & \text{if $n\equiv 2 ~ (mod\, 3)$ and $i=1$ or $\frac{n}{2}$ or $n-1$} \\
s(P_n)+1, & \text{if $n\equiv 1 ~ (mod\, 3)$ and $i=2$ or $n-2$} \\
s(P_n), & \text{otherwise.}\\
\end{array} \right.
\end{equation*}
    \end{enumerate}
\end{proposition}

\proof
 \begin{enumerate}
      \item [(i)] It is clear that $C_n-e=P_n$. Considering the saturation number of $P_n$ and $C_n$, we have the result.
       \item[(ii)] For every $e\in E(P_n)$, we have ${P_n}-e={P_m}\cup {P_t}$ where $m,t\in \Bbb{N}$ and $m+t=n$. By Lemma \ref{union}, we can conclude that $s({P_n}-e)=s(P_m)+s(P_t)$. Now we cosider three following cases: \\
Case 1) If $n=3k$ for some $k\in \Bbb{Z}$, then $m=3k_1$ and $t=3k_2$ where $k_1, k_2\in \Bbb{Z}$ and $k_1+k_2=k$.  So $s({P_n}-e)=s(P_{3k_1})+s(P_{3k_2})=k_1+k_2=k=s(P_n)$; or $m=3k_1+1$ and $t=3k_2+2$ where $k_1, k_2\in \Bbb{Z}$ and $k_1+k_2=k-1$. Then $s({P_n}-e)=s(P_{3k_1+1})+s(P_{3k_2+2})=k_1+(k_2+1)=k=s(P_n).$ \\
Case 2) Suppose that $n=3k+1$ for some $k\in \Bbb{Z}$. So we have $m=3k_1$ and $t=3k_2+1$ where $k_1, k_2\in \Bbb{Z}$ and $k_1+k_2=k$. Therefore $s({P_n}-e)=s(P_{3k_1})+s(P_{3k_2+1})=k_1+k_2=k=s(P_n)$; or $m=3k_1+2$ and $t=3k_2+2$ where $k_1, k_2\in \Bbb{Z}$ and $k_1+k_2=k-1$. Then $s({P_n}-e)=s(P_{3k_1+2})+s(P_{3k_2+2})=(k_1+1)+(k_2+1)=k+1=s(P_n)+1.$ Note that it case happens when $e=e_2$ or $e=e_{n-2}$. \\
Case 3) If $n=3k+2$ for some $k\in \Bbb{Z}$, then $m=3k_1+1$ and $t=3k_2+1$ where $k_1, k_2\in \Bbb{Z}$ and $k_1+k_2=k$ and we have $s({P_n}-e)=s(P_{3k_1+1})+s(P_{3k_2+1})=k_1+k_2=k=s(P_n)-1.$
Note that this case happens if $e=e_1$, $e=e_{n-1}$ or $e=\frac{n}{2}$ (when $n$ is even); or $m=3k_1$ and $t=3k_2+2$ where $k_1, k_2\in \Bbb{Z}$ and $k_1+k_2=k$. Then $s({P_n}-e)=s(P_{3k_1})+s(P_{3k_2+2})=
k_1+(k_2+1)=k+1=s(P_n).$
\qed
\end{enumerate}

\subsection{Saturation number of the corona product of two specific graphs} 
  
In this subsection, we consider the corona product of two certain graphs and study their saturation number. We recall that the corona of two graphs $G_1$ and $G_2$, written as ${G_1}\circ{G_2}$, is the graph obtained by taking one copy of $G_1$ and $|V(G_1)|$ copies of $G_2$, and then joining the $i$th vertex of $G_1$ to every vertex in the $i$th copy of $G_2$. 

\begin{theorem}\label{corona}
Let $G$ be a simple graph of order $n$. Then 
$$s(G\circ \overline{K_m})=\alpha^{\prime}(G)+l$$ 
where $\alpha^{\prime}(G)$ is the size of maximum matching $M$ of  $G$ and $l$ is the number of $M$-unsaturated vertices of $G$. In addition, if $G$ has a perfect matching, then 
$$s(G\circ \overline{K_m})=\frac{n}{2}.$$
\end{theorem}

\proof 
Suppose that $M$ is a maximum matching of $G$ and $S$ has the smallest cardinality over all maximal matchings of 
$G\circ \overline{K_m}$. Let the vertex $u$ be a $M$-unsaturated vertex of $G$ and $(\overline{K_m})_u$ be a copy of 
$\overline{K_m}$  with the vertex set $\{v_1,...,v_m\}$ such that the vertex $u$ is adjacent to all vertices of 
$(\overline{K_m})_u$. Thus there exist $v\in \{v_1,...,v_m\}$ such that $uv\in S$. Now if for $M$-unsaturated vertices $u_1, u_2, ..., u_l$ of $G$, put 
$$S=M\cup\big\{u_iv|1\leq i\leq l, v\in V\big({(\overline{K_m})}_{u_i}\big)\big\},$$ then $S$ has the smallest cardinality over all maximal matchings of $G\circ \overline{K_m}$ and so we have the result. \\
Now suppose that the graph $G$ has a perfect matching $N=\{e_1, e_2, ..., e_{\frac{n}{2}}\}$. To obtain a maximal matching $S$ with the smallest cardinality of $G\circ \overline{K_m}$, it is enough to put $S=N$, and so 
$s(G\circ \overline{K_m})=|S|=|N|=\frac{n}{2}$. 
\qed

\begin{corollary} 
For every $n,m\geq 3$, $s({P_n}\circ \overline{K_m})=s({C_n}\circ \overline{K_m})=\left\lceil{\frac{n}{2}}\right\rceil$. 
\end{corollary} 

\proof 
Clearly, for every even $n$, the path $P_n$ and the cycle $C_n$ have perfect matching and by Theorem \ref{corona} we have $s({P_n}\circ \overline{K_m})=s({C_n}\circ \overline{K_m})=\frac{n}{2}$. Now suppose that $n$ is odd. Then for every maximum matching $M$, the path $P_n$ and the cycle $C_n$ have only one $M$-unsaturated vertex. So by 
Theorem \ref{corona}, $$s({P_n}\circ \overline{K_m})=s({C_n}\circ \overline{K_m})=\left\lfloor{\frac{n}{2}}
\right\rfloor+1=\left\lceil{\frac{n}{2}}\right\rceil.$$    
\qed

\begin{theorem}\label{GPm}
For any simple graph $G$ of order $n$, we have $$s(G\circ {P_m})=\left\{
\begin{array}{lcl}
ns(P_m)+\alpha^{\prime}(G)+l, & \quad\mbox{if $m\equiv 1 ~ (mod\, 3)$}\\[15pt]
ns(P_m), & \quad\mbox{otherwise.}
\end{array}
\right.$$ where $\alpha^{\prime}(G)$ is the size of maximum matching $M$ of $G$ and $l$ is the number of 
$M$-unsaturated vertices of $G$.
\end{theorem}

\proof
Let $u\in V(G)$ and ${(P_m)}_u$ be a copy of $P_m$ with the vertex set $\{v_1,...,v_m\}$ such that the vertex $u$ is adjacent to all vertices of ${(P_m)}_u$. Suppose that  $S$ and $S_1$ have the smallest cardinality over all maximal matchings of graphs $G\circ P_m$ and $P_m$, respectively. So $s(G\circ P_m)=|S|$ and $s(P_m)=|S_1|$. We have three following cases: 

\medskip
Case 1) Suppose that $m=3k$ for some $k\in \Bbb{Z}$. It can easily verified that if $v\in \{v_2, v_{m-1}\}$, then the 
edge $uv$ belongs to $S$. Hence every vertex of $G$ is $S$-saturated and so there is no edge of $G$ in $S$. 
Also the number of edges in the smallest maximal matching of the graph ${(P_m)}_u-v$ with the edge 
$uv$ is equal to the number of edges in the smallest maximal matching of the graph $P_m$. Thus we can conclude that 
$|S|=n|S_1|$ and we have $$s(G\circ P_m)=|S|=n|S_1|=ns(P_m).$$
 
Case 2) Assume that $m=3k+1$ for some $k\in \Bbb{Z}$. If $M$ is a perfect matching of $G$, then put $S=\underbrace{S_1\cup S_1\cup ...\cup S_1}_{n-times}\cup M$. So $$s(G\circ P_m)=|S|=n|S_1|+|M|=ns(P_m)+\alpha^{\prime}(G).$$ 
Now suppose that $u$ is an $M$-unsaturated vertex of $G$ and ${(P_m)}_u$ is a copy of $P_m$ with the vertex set 
$\{v_1,...,v_m\}$ such that the vertex $u$ is adjacent to all vertices of ${(P_m)}_u$. Then there exist 
$v\in \{v_1,...,v_m\}$ such that $uv\in S$. Suppose that $u_1, u_2, ..., u_l$ are $M$-unsaturated vertices of $G$ and put  
$$S=\underbrace{S_1\cup S_1\cup ...\cup S_1}_{n-times}\cup M\cup \big\{u_iv|1\leq i\leq l, v\in V\big({(P_m)}_{u_i}\big)\big\}.$$ Then we have $$s(G\circ P_m)=|S|=n|S_1|+|M|+l=ns(P_m)+\alpha^{\prime}(G)+l.$$ 

Case 3) If $m=3k+2$ for some $k\in \Bbb{Z}$ and $v$ is one of the vertices in the set $\{v_1, v_2, v_4, v_5, v_7, ..., v_{m-3}, v_{m-1}, v_{m}\}$, then the edge $uv$ belongs to $S$.
 Therefore, every vertex of $G$ is $S$-saturated and so there is no edge of $G$ in $S$.
Similar to  proof of Case 1, since the number of edges in the smallest maximal matching of the graph 
${(P_m)}_u-v$ with the edge $uv$ is equal to the number of edges in the smallest maximal matching of the path  
$P_m$, so $|S|=n|S_1|$ and it implies that $$s(G\circ P_m)=|S|=n|S_1|=ns(P_m).$$ 
\qed

\begin{theorem}\label{GCm}
For any simple graph  $G$ of order $n$, we have
$$s(G\circ {C_m})=\left\{
\begin{array}{lcl}
ns(C_m)+\alpha^{\prime}(G)+l, & \quad\mbox{if  $m\equiv 0 ~ (mod\, 3)$}\\[15pt]
ns(C_m), & \quad\mbox{otherwise.}
\end{array}
\right.$$ where $\alpha^{\prime}(G)$ is the size of maximum matching $M$ of the graph $G$ and $l$ is the number of 
$M$-unsaturated vertices of $G$.
\end{theorem}

\proof
Let $u\in V(G)$ and $(C_m)_u$ be a copy of $C_m$ such that the vertex $u$ is adjacent to every vertex of $(C_m)_u$. 
Suppose that $S$ and $S_1$ have the smallest cardinality over all maximal matchings of graphs $G\circ C_m$ and 
$C_m$, respectively. Then $s(G\circ C_m)=|S|$ and $s(C_m)=|S_1|$. 
If $m=3k+1$ or $m=3k+2$ for some $k\in \Bbb{Z}$, then the edge $uv$ belongs to $S$, where $v$ is a vertex of $C_m$. Thus every vertex of $G$ is $S$-saturated and there is no edge of $G$ in $S$. Also we have 
$s(C_m-v)=s(P_{m-1})=|S_1|-1$. Then we can conclude that $$s(G\circ {C_m})=|S|=n|S_1|=ns(C_m).$$ 

Now suppose that $m=3k$ for some $k\in \Bbb{Z}$. Similar to the proof of Theorem \ref{GPm}, if $M$ is a perfect matching of $G$, then put $S=\underbrace{S_1\cup S_1\cup ...\cup S_1}_{n-times}\cup M$ and so we have the result. But if $G$ does not have a perfect matching, then assume that $u$ is an $M$-unsaturated vertex of $G$ and ${(C_m)}_u$ is a copy of $C_m$ with the vertex set $\{v_1, ..., v_m\}$ such that the vertex $u$ is adjacent to all vertices of ${(C_m)}_u$. Then there exist $v\in \{v_1, ..., v_m\}$ such that $uv\in S$. Now for $M$-unsaturated vertices $u_1, u_2, ..., u_l$ of $G$, it is enough to put $$S=\underbrace{S_1\cup S_1\cup ...\cup S_1}_{n-times}\cup M\cup \big\{u_iv| 1\leq i\leq l,  v\in V\big({(C_m)}_{u_i}\big)\big\}$$ and so the result follows.
\qed 


\begin{proposition}\label{K_1G}
For any simple graph $G$, $s(G)\leq s(K_1\circ G)\leq 1+s(G).$
\end{proposition}

\proof 
Suppose that $K_1=\{u\}$ and $S_2$ has the smallest cardinality over all maximal matchings of graph 
$G-v$ where $v\in V(G)$. So $s(G-v)=|S_2|$.
To obtain the maximal matching $S$ with the smallest cardinality of $K_1\circ G$, put $S=\{uv\}\cup S_2$. Then we have 
$$s(K_1\circ G)=|S|=1+|S_2|=1+s(G-v).$$ Also it can easily verified that $s(G-v)=s(G)$ or 
$s(G-v)=s(G)-1$ and it implies the result.
\qed  

\medskip
As some examples, we state the saturation number of $K_1\circ P_n$, $K_1\circ C_n$ and $K_1\circ W_n$ in the following.

\begin{example}  
	For every $n\in \Bbb{N}$, we have
    \begin{enumerate}
           \item [(i)] $s(K_1\circ P_n)=\left\{
\begin{array}{lcl}
1+s(P_n), & \quad\mbox{if  $n\equiv 1 ~ (mod\, 3)$}\\[15pt]
s(P_n), & \quad\mbox{otherwise.}
\end{array}
\right.$
           \item [(ii)] $s(K_1\circ C_n)=\left\{
\begin{array}{lcl}
1+s(C_n), & \quad\mbox{if  $n\equiv 0 ~ (mod\, 3)$}\\[15pt]
s(C_n), & \quad\mbox{otherwise.}
\end{array}
\right.$
           \item [(iii)] $s(K_1\circ W_n)=\left\{
\begin{array}{lcl}
1+s(W_n), & \quad\mbox{if $n\equiv 0 ~ (mod\, 3)$}\\[15pt]
s(W_n), & \quad\mbox{otherwise.}
\end{array}
\right.$
    \end{enumerate}
    \end{example} 
\medskip

\begin{proposition}
For any simple graph $G$, $s(\overline{K_m}\circ G)=m.s(K_1\circ G).$  
\end{proposition}

\proof
It is clear that $\overline{K_m}\circ G=\underbrace{(K_1\circ G)\cup ...\cup (K_1\circ G)}_{m-times}$. Now the result follows from Lemma \ref{union}.
\qed  

\medskip
We end this section with the following corollary which is an immediate consequence of Theorems \ref{corona}, \ref{GPm} and \ref{GCm}.   

\begin{corollary}
Let $G_1$ and $G_2$ be two simple graphs and $|V(G_1)|=n$. 
Then  $$ns(G_2)\leq s(G_1\circ G_2)\leq ns(G_2)+\alpha^{\prime}(G_1)+l$$
where $\alpha^{\prime}(G_1)$ is the size of maximum matching $M$ of the graph $G_1$ and $l$ is the number of 
$M$-unsaturated vertices of $G_1$. 
\end{corollary}

\section{Saturation number of the link and the chain of graphs} 

\begin{figure}[h]
	\begin{center}
		\begin{minipage}{5cm}
			\hspace{1cm}
			\includegraphics[width=5cm,height=.9cm]{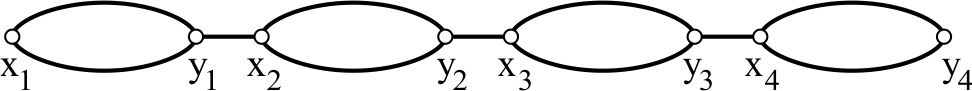}
		\end{minipage}
		\begin{minipage}{5cm}
			\includegraphics[width=5cm,height=1cm]{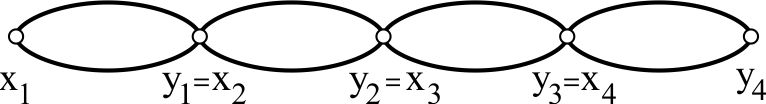}
		\end{minipage}
		\caption{ \label{1} \small A link and a chain of four graphs, respectively. }
	\end{center}
\end{figure}

In this section we investigate the saturation number of some graphs with specific construction that are of importance in chemistry. First we consider the link of graphs. Let $G_1$, $G_2$, ..., $G_k$ be a finite sequence of pairwise disjoint connected graphs and let $x_i, y_i \in V(G_i)$. The link $G$ of the graphs $\{G_i\}_{i=1}^k$ with respect to the vertices 
$\{x_i, y_i\}_{i=1}^k$ is obtained by joining by an edge the vertex $y_i$ of $G_i$ with the vertex $x_{i+1}$ of $G_{i+1}$ for all $i=1,2,...,k-1$ {\rm(}see Figure \ref{1} for $k=4${\rm)} \cite{Deu}.

In general attitude, there is no a certain relation between $s\big(L(G_1, G_2, ..., G_k)\big)$ and $s(G_i)$ where 
$1\leq i\leq k$. For example, $s\big(L(C_4, P_5, P_4)\big)=s(C_4)+s(P_5)+s(P_4)$, $s\big(L(C_4, P_5, C_4)\big)=2s(C_4)+s(P_5)-1$ and $s\big(L(C_4, P_5, C_4, P_5)\big)=2s(C_4)+2s(P_5)-2$.
But for some special graphs there are certain relations. In the following theorem, we consider the saturation number of the link of some well-known graphs. As usual $d(x,y)$ denotes the distance between two vertices $x$ and $y$. 

\begin{proposition}\label{lpc}
    \begin{enumerate}
           \item [(i)] Let $P_m$ be a path of order $m$. Then 
$$s\big(L(\underbrace{P_m, P_m, ..., P_m}_{n-times})\big)=\left\{
\begin{array}{lcl}
ns(P_m), & \quad\mbox{if $m\equiv 0 ~ (mod\, 3)$}\\[15pt]
ns(P_m)+\left\lceil{\frac{n-1}{3}}\right\rceil, & \quad\mbox{if $m\equiv 1 ~ (mod\, 3)$}\\[15pt]
ns(P_m)-\left\lceil{\frac{n-1}{3}}\right\rceil, & \quad\mbox{if $m\equiv 2 ~ (mod\, 3)$}.
\end{array}
\right.$$

           \item [(ii)] Let $C_m$ be a cycle of order $m$ and $x_i, y_i\in V({(C_m)}_i)$ for every $1\leq i\leq n$. If $d(x_i,y_i)=d_i$, then we have $$s\big(L(\underbrace{C_m,..., C_m}_{n-times})\big)=\left\{
\begin{array}{lcl}
ns(C_m), & \quad\mbox{if $m\equiv 0 ~ (mod\, 3)$, $1\leq d_i\leq 5$ }\\[15pt]
ns(C_m)-\left\lfloor{\frac{n}{2}}\right\rfloor, & \quad\mbox{if $m\equiv 1 ~ (mod\, 3)$, $d_i\in \{1,3.4\}$ }\\[15pt]
ns(C_m)-(n-1), & \quad\mbox{if $m\equiv 1 ~ (mod\, 3)$, $d_i\in \{2,5\}$ }\\[15pt]
ns(C_m)-\left\lceil{\frac{n-1}{3}}\right\rceil, & \quad\mbox{if $m\equiv 2 ~ (mod\, 3)$, $d_i\in \{1,4\}$ }\\[15pt]
ns(C_m)-\left\lfloor{\frac{n}{2}}\right\rfloor, & \quad\mbox{if $m\equiv 2 ~ (mod\, 3)$, $d_i\in \{2,3,5\}$}.
\end{array}
\right.$$
\end{enumerate}
\end{proposition} 

\proof
    \begin{enumerate}
           \item [(i)] Clearly, the link of $n$ paths of order $m$ is a path with $mn$ vertices. Suppose that $\{v_1,..., v_m\}$ is the vertex set of $P_m$. It is enough to choose $S=\{v_2v_3, v_5v_6, v_8v_9, ...\}$. Then $S$ has the smallest cardinality over all maximal matchings of $P_m$. Now the result can be easily verified. 

           \item [(ii)] We consider all cases which has stated for the saturation number of the link of some cycles. The construction of maximal matchings with the smallest cardinality for examples of  these cases, has shown in Figures 
2, 3, 4, 5 and 6, respectively. So we have the result. \qed 

\begin{figure}[ht]
\centerline{\includegraphics[width=8cm]{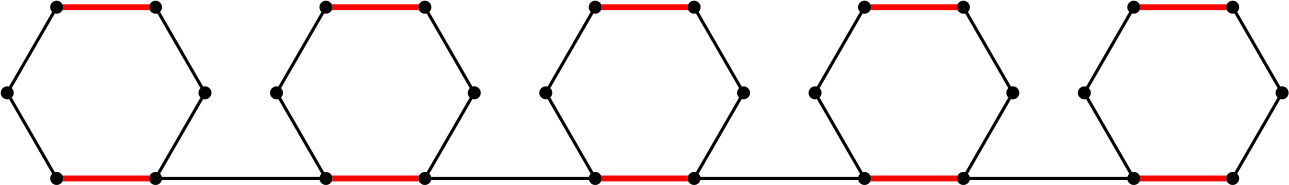}}
\medskip
\centerline{\includegraphics[width=8cm]{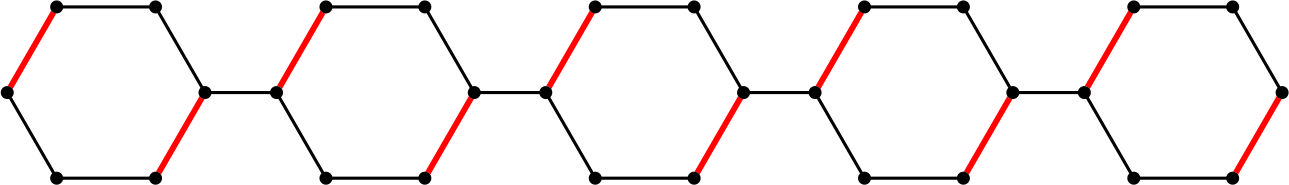}}
\medskip
\centerline{\includegraphics[width=8cm]{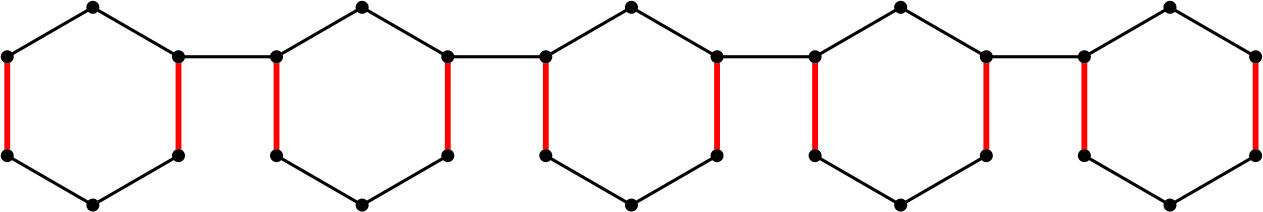}}
\caption{ \label{2} \small $s\big(L(C_6, C_6, C_6, C_6, C_6)\big)=10$.}
\end{figure}

\begin{figure}[ht]
\centerline{\includegraphics[width=8cm]{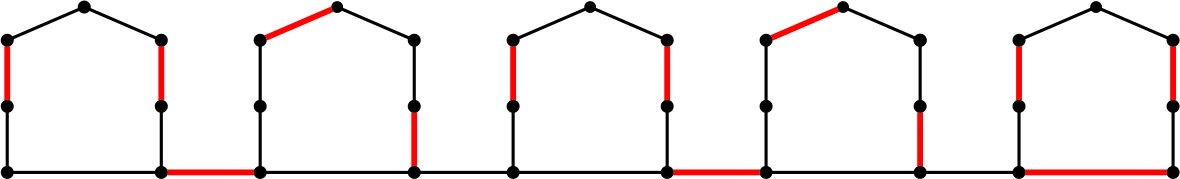}}
\medskip
\centerline{\includegraphics[width=8cm]{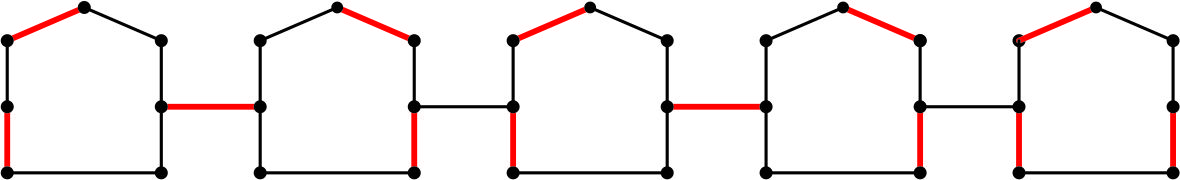}}
\caption{ \label{2} \small $s\big(L(C_7, C_7, C_7, C_7, C_7)\big)=13$.}
\end{figure}

\begin{figure}[ht]
\centerline{\includegraphics[width=8cm]{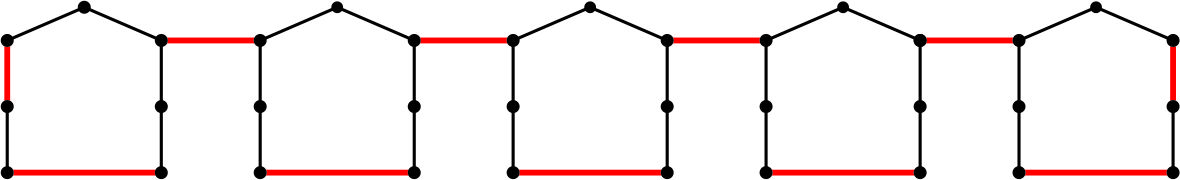}}
\caption{\label{4}\small $s\big(L(C_7, C_7, C_7, C_7, C_7)\big)$=11.}
\end{figure}

\begin{figure}[h]
\centerline{\includegraphics[width=8cm]{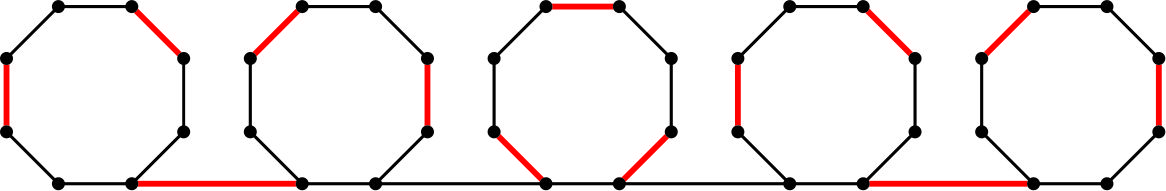}}
\medskip
\centerline{\includegraphics[width=8cm]{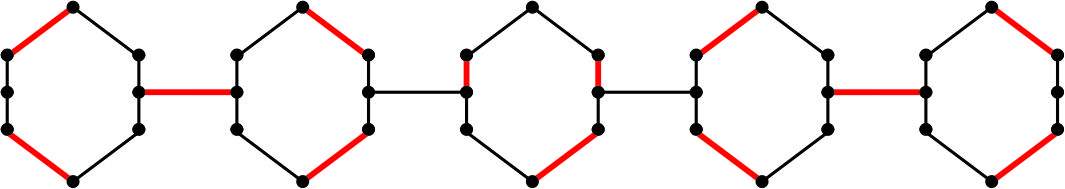}}
\caption{\label{4}\small $s\big(L(C_8, C_8, C_8, C_8, C_8)\big)$=13.}
\end{figure}

\begin{figure}[h]
\centerline{\includegraphics[width=8cm]{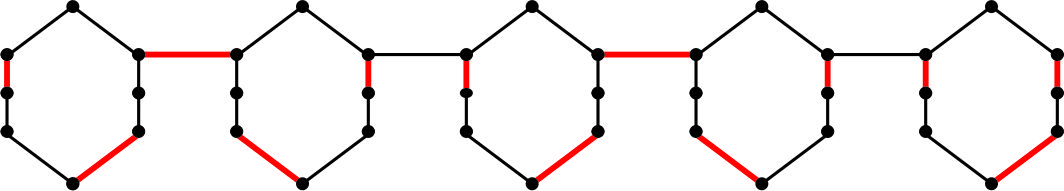}}
\medskip
\centerline{\includegraphics[width=8cm]{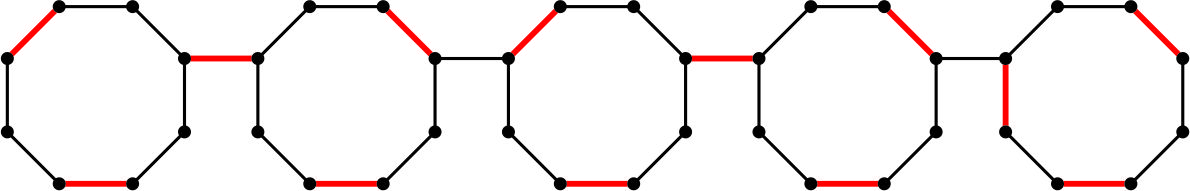}}
\caption{\label{4}\small $s\big(L(C_8, C_8, C_8, C_8, C_8)\big)$=13.}
\end{figure}
\end{enumerate}

Here we study the saturation number of the chain of graphs. Let $G_1, G_2, ..., G_k$ be a finite sequence of pairwise disjoint connected graphs and let $x_i,y_i\in V(G_i)$. The chain $G$ of the graphs $\{G_i\}_{i=1}^{k}$ with respect to the vertices $\{x_i,y_i\}_{i=1}^{k}$ is obtained by identifying the vertex $y_i$ with the vertex $x_{i+1}$ for $1\leq i\leq k-1$ {\rm(}see Figure \ref{1} for $k =4${\rm)} \cite{Deu}. Similar to Proposition \ref{lpc}, the following observation presents  the saturation number of the chain of paths and cycles, respectively.

\begin{observation}\label{chain}
    \begin{enumerate}
           \item [(i)] Let $P_m$ be a path of order $m$.  Then 
$$s\big(C(\underbrace{P_m, P_m, ..., P_m}_{n-times})\big)=\left\{
\begin{array}{lcl}
ns(P_m)-\left\lfloor{\frac{n}{3}}\right\rfloor, & \quad\mbox{if  $m\equiv 0 ~ (mod\, 3)$}\\[15pt]
ns(P_m), & \quad\mbox{if    $m\equiv 1 ~ (mod\, 3)$}\\[15pt]
ns(P_m)+\left\lfloor{\frac{-2(n-1)}{3}}\right\rfloor, & \quad\mbox{if  $m\equiv 2 ~ (mod\, 3)$}.
\end{array}
\right.$$
           \item [(ii)] Let $C_m$ be a cycle of order $m$ and $x_i, y_i\in V({(C_m)}_i)$ for every $1\leq i\leq n$. If $d(x_i,y_i)=d_i$, then we have $$s\big(C(\underbrace{C_m,..., C_m}_{n-times})\big)=\left\{
\begin{array}{lcl}
ns(C_m)-\left\lfloor{\frac{n-1}{2}}\right\rfloor, & \quad\mbox{if $m\equiv 0 ~ (mod\, 3)$, $d_i\in \{1,2,4,5\}$ }\\[15pt]
ns(C_m), & \quad\mbox{if $m\equiv 0 ~ (mod\, 3)$, $d_i=3$ }\\[15pt]
ns(C_m)-(n-1), & \quad\mbox{if $m\equiv 1 ~ (mod\, 3)$, $1\leq d_i\leq 5$ }\\[15pt]
ns(C_m)-\left\lceil{\frac{n}{2}}\right\rceil, & \quad\mbox{if $m\equiv 2 ~ (mod\, 3)$, $d_i\in \{1,4\}$ }\\[15pt]
ns(C_m)-(n-1), & \quad\mbox{if $m\equiv 2 ~ (mod\, 3)$, $d_i\in \{2,3,5\}$}.
\end{array}
\right.$$
    \end{enumerate}
\end{observation}

\begin{figure}[h]
	\begin{center}
		\begin{minipage}{5cm}
			\hspace{4cm}
			\includegraphics[width=5cm,height=1.5cm]{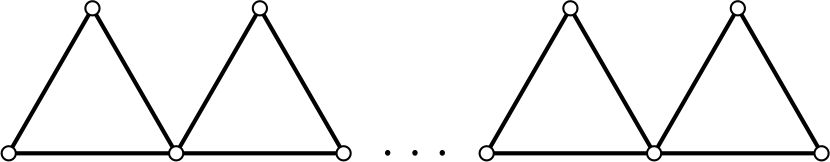}
		\end{minipage}
		\begin{minipage}{5cm}
			\includegraphics[width=5cm,height=2cm]{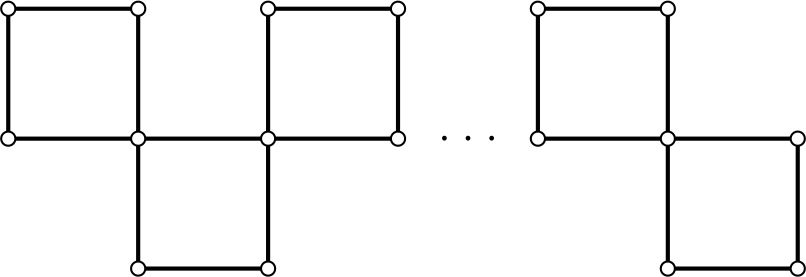}
		\end{minipage}
		\caption{ \label{8} \small A chain triangular cactus $T_n$ and  square cactus $O_n$,  respectively.}
	\end{center}
\end{figure}

Now, we shall study the saturation number of cactus graphs that are of importance in chemistry. 
A cactus graph is a connected graph in which no edge lies in more than one cycle. Consequently, each block of a cactus graph is either an edge or a cycle. If all blocks of a cactus $G$ are cycles of the same size $k$, the cactus is $k$-uniform. 
A triangular cactus is a graph whose blocks are triangles, i.e., a $3$-uniform cactus. A vertex shared by two or more triangles is called a cut-vertex. If each triangle of a triangular cactus $G$ has at most two cut-vertices, and each cut-vertex is shared by exactly two triangles, we say that $G$ is a chain triangular cactus. The number of triangles in $G$ is called the length of the chain. An example of a chain triangular cactus is shown in Figure \ref{8}. 

Obviously, all chain triangular cactus of the same length are isomorphic. Hence, we denote the chain triangular cactus of length $n$ by $T_n$. clearly, a chain triangular cactus of length $n$ has $2n+1$ vertices and $3n$ edges \cite{Jah}. 
By replacing triangles in chain triangular chain $T_n$ by cycles of length $4$, we obtain cactus whose every block is $C_4$ as shown in Figure \ref{8}. We call such cactus, square cactus and denote a chain square cactus of length $n$ by $O_n$ \cite{Jah}. The following theorem gives the saturation number of the chain triangular cactus and a chain square cactus of length $n$.

\begin{theorem}
Let $T_n$ and $O_n$ be a chain triangular cactus and a chain square cactus of length $n$, respectively. Then we have 
     \begin{enumerate}
           \item [(i)] For every $n\geq 1$, $s(T_n)=\left\lfloor\frac{n-2}{2}\right\rfloor+2$.
           \item [(ii)] For every $n\geq 1$, $s(O_n)=n+1$.
     \end{enumerate}
\end{theorem}

\proof
     \begin{enumerate}
           \item [(i)] By the Observation \ref{chain} (ii), $s(T_n)=s\big(C(\underbrace{C_3, C_3, ..., C_3}_{n-times})\big)=n-\left\lfloor{\frac{n-1}{2}}\right\rfloor$. Now the induction on $n$ implies that 
$n-\left\lfloor{\frac{n-1}{2}}\right\rfloor=\left\lfloor{\frac{n-2}{2}}\right\rfloor+2$ and so we have the result.

           \item [(ii)] It follows from Observation \ref{chain} (ii).
     \end{enumerate}
\qed



\begin{thebibliography}{1}

      \bibitem{Jah} S. Alikhani, S. Jahari, M. Mehryar and R. Hasni, {\it Counting the number of dominating sets of cactus chains}, Optoelectron. Adv. Mater.- Rapid Comm, 8, no. 9-10 (2014) 955-960.

      \bibitem{Ves} V. Andova, F. Kardo\v s and R. \v Skrekovski, {\it Sandwiching saturation number of fullerene graphs}, arXiv:1405.2197 (2014).


      \bibitem{Cyv} S. J. Cyvin and I. Gutman, {\it Kekul\'e structures in benzenoid hydrocarbons}, volume 46 of Lecture
Notes in Chemistry, Springer Science, Heidelberg, 1988.

     \bibitem{Deu} E.  Deutsch and S. Klav\v zar, {\it Computing the Hosoya Polynomial of Graphs from Primary Subgraphs}, MATCH Commun. Math. Comput. Chem. 70 (2013) 627-644.

      \bibitem{Ful} T. Do\v sli\'c, {\it On lower bounds of number of perfect matchings in fullerene graphs}, J. Math. Chem.
24 (1998), 359-364.

     \bibitem{Sli} T. Do\v sli\'c, {\it Fullerene graphs with exponentially many perfect matchings}, J. Math. Chem. 41
(2007), 183-192.

    \bibitem{Sat} T. Do\v sli\'c, {\it Saturation number of fullerene graphs}, J. Math. Chem. 43 (2008), 647-657.

     \bibitem{Dos}  T. Do\v sli\'c and I. Zubac, {\it Saturation number of benzenoid graphs}, MATCH Commun. Math. Comput. Chem 73 (2015) 491-500.

      \bibitem{Edm} J. Edmonds, {\it Paths, trees, and flowers}, Canad. J. Math. 17 (1965) 449-467.

     \bibitem{Fad} J. Faudree, R. J. Faudree, R. J. Gould and M. S. Jacobson, {\it Saturation numbers for trees}, Electron.
J. Combin. 16 (2009).

     \bibitem{Fre} A. Frendrup, B. Hartnell and P. D. Vestergaard, {\it A note on equimatchable graphs}, Australas. J. Combin. 46 (2010), 185-190.

     \bibitem{Har} F. Harary, {\it Graph Theory}. Reading, MA: Addison-Wesley, 1994.   

     \bibitem{Kar} F. Kardo\v s, D. Kr\'al, J. Mi\v skuf and J. S. Sereni, {\it  Fullerene graphs have exponentially many perfect
matchings}, J. Math. Chem. 46 (2009), 443-447.

     \bibitem{Kla} M. Klazar, {\it Twelve countings with rooted plane trees}, European J. Combin. 18 (1997), 195-210.

     \bibitem{Lov} L. Lov\'asz and M.D. Plummer, {\it Matching Theory}, Annals of Discrete Math. Vol. 29, North-Holland, Amsterdam, 1986.

     \bibitem{Wan} S. G. Wagner, {\it On the number of matchings of a tree}, European J. Combin. 28 (2007), 1322-1330.

     \bibitem{Wes} D. B. West, {\it Introduction to Graph Theory}, Prentice Hall, Inc., Upper Saddle River, NJ, 1996.

     \bibitem{Zha} H. Zhang and F. Zhang, {\it New lower bound on the number of perfect matchings in fullerene
graphs}, J. Math. Chem. 30 (2001), 343-347.

\end{thebibliography}
\end{document}